\documentstyle[12pt]{article}
\begin{document}
\newtheorem{Def}{Definition}
\newtheorem{thm}{Theorem}
\newtheorem{lem}{Lemma}
\newtheorem{rem}{Remark}
\newtheorem{prop}{Proposition}
\newtheorem{cor}{Corollary}
\title
{A Liouville type theorem
 for some conformally invariant fully nonlinear equations } 
\author{\ Aobing Li
\ \ \ \ \& \ \ \ YanYan Li\thanks{Partially supported by
NSF Grant DMS-0100819.}
\\ Department of Mathematics\\ Rutgers University\\
110 Frelinghuysen Rd.\\
Piscataway, NJ 08854
}
\date{}
\maketitle
\input { amssym.def}

Following the approach in our earlier paper \cite{LL}
and using the gradient estimates developed in
\cite{LL} and \cite{LL2}, we
give another  Liouville type theorem 
 for some conformally invariant fully nonlinear equations.
Various Liouville type theorems for conformally invariant equations
have been obtained by Obata,  Gidas-Ni-Nirenberg,
Caffarelli-Gidas-Spruck, Viaclovsky, Chang-Gursky-Yang, and Li-Li.  For
these, as well as for related works, see
\cite{LL} and the references therein.

For $n\ge 3$, 
let ${\cal S}^{n\times n}$ be the set of  $n\times n$ real symmetric
matrices,   
${\cal S}^{n\times n}_+\subset {\cal S}^{n\times n}$ be the set
of positive definite matrices, and let $O(n)$ be 
the set of $n\times n$ real orthogonal matrices.

For  $1\le k\le n$,
let
$$
\sigma_k(\lambda)=\sum_{1\le i_1<\cdots<i_k\le n}\lambda_{i_1}
\cdots \lambda_{i_k},\qquad\qquad \lambda=(\lambda_1,
 \cdots, \lambda_n)\in \Bbb R^n,
$$
denote the $k-$th symmetric function, and let
 $\Gamma_k$ denote the
connected component of $\{\lambda\in \Bbb R^n\ |\
\sigma_k(\lambda)>0\}$ containing  the
positive cone $\{\lambda\in \Bbb R^n\ |\
\lambda_1, \cdots, \lambda_n>0\}$.
It is  known that
$$
\Gamma_n=\{\lambda\in \Bbb R^n\ |\
\lambda_1, \cdots, \lambda_n>0\}, \qquad
\Gamma_1=\{\lambda\in \Bbb R^n\ |\
\lambda_1+ \cdots + \lambda_n>0\},
$$
$$
\Gamma_k=\{\lambda\in \Bbb R^n\ |\ \sigma_1(\lambda)>0, \cdots,
\sigma_k(\lambda)>0\},
$$
$\Gamma_k$ is
a convex cone with its vertex at the origin with the properties
$$
\Gamma_n\subset \cdots\subset \Gamma_2\subset \Gamma_1,
$$
$$
\frac{\partial \sigma_k}{\partial \lambda_i}>0\quad
\mbox{in}\ \Gamma_k, \ 1\le i\le n,
$$
$$
\sigma_k^{\frac 1k}\ \mbox{is concave in}\ \Gamma_k.
$$

For a positive $C^2$ function $u$, let
$$
A^u:= -\frac{2}{n-2}u^{  -\frac {n+2}{n-2} }
\nabla^2u+ \frac{2n}{(n-2)^2}u^ { -\frac {2n}{n-2} }
\nabla u\otimes\nabla u-\frac{2}{(n-2)^2} u^ { -\frac {2n}{n-2} }
|\nabla u|^2I,
$$
where  $I$ is the $n\times n$ identity matrix.

Assume $U\subset {\cal S}^{n\times n}$  is an open set satisfying
\begin{equation}
O^{-1}UO=U,\qquad\forall\ O\in O(n),
\label{U1}
\end{equation}
and
\begin{equation}
U\cap \{M+tN\ |\ 0<t<\infty\}\ \mbox{is convex}\qquad \forall\
M\in {\cal S}^{n\times n}, N\in {\cal S}^{n\times n}_+, 
\label{U2}  
\end{equation}
and
\begin{equation}
\Gamma_U:=\{\lambda(M)|~M\in U\}\subset\Gamma_k,\quad\mbox{for
some}~k>\frac{n+1}{2},
\label{Gammak}
\end{equation}
where $\lambda(M)$ denotes the eigenvalues of $M$.

Let $F\in C^2(U)$ satisfy
\begin{equation}
F(O^{-1}MO)=F(M),\qquad \forall\
M\in U, \  O\in O(n),
\label{F1}
\end{equation}
\begin{equation}
0\ \mbox{does not belong to }\ F^{-1}(1),
\label{F3}
\end{equation}
\begin{equation}
\left(F_{ij}(M)\right)>0,\qquad \forall\ M\in U,
\label{F2}
\end{equation}
\begin{equation}
F\ \mbox{is locally concave in}\ U,
\label{F4}
\end{equation}
and, for some $0<\gamma\le 1$,
\begin{equation}
\sum_{i.j=1}^n F_{ij}(M)M_{ij}\le \frac 1\gamma|M|^{1-\gamma}\sum_{i=1}^n
F_{ii}(M), 
\quad
\forall\ M\in U, F(M)=1, |M|\ge 1,
\label{F5}
\end{equation}
where $F_{ij}(M):=\frac{\partial F}{ \partial M_{ij} }(M)$.

We establish in this paper the following Liouville type theorem.
\begin{thm}  For $n\ge 3$, let $U\subset {\cal S}^{n\times n}$ be an 
open set satisfying (\ref{U1}), (\ref{U2})
and (\ref{Gammak}), and let
$F\in C^2(U)$  satisfy  (\ref{F1}), (\ref{F3}),
(\ref{F2}), (\ref{F4}) and (\ref{F5}). 
Let $u\in C^4(\Bbb R^n)$ be  a positive solution  of
\[
F(A^u)=1,\quad A^u\in U,\quad\mbox{on}~\Bbb{R}^n.
\]
Then for some $\bar{x}\in\Bbb{R}^n$,
 and  some positive constants $a$ and $b$
satisfying
$2b^2a^{-2}I\in U$ and
$F(2b^2a^{-2}I)=1$,
\begin{equation}
u(x)\equiv \left(\frac a {1+b^2|x-\bar{x}|^2}\right)^{\frac{n-2}{2}},
\qquad \forall\ x\in \Bbb R^n.
\label{form}
\end{equation}
\label{highk}
\end{thm} 

\begin{rem} In Theorem \ref{highk}, if $F$ is in $C^{2,\beta}(U)$
for some $\beta\in (0,1)$, then, since the equation is elliptic,
any positive $C^2$ solution $u$ is in fact in $C^{4,\beta}$.
\end{rem}

We give a consequence of Theorem \ref{highk}.

Let
\begin{equation}
\Gamma\subset \Bbb R^n \ \mbox{be
an open convex cone  with its vertex at the origin }
\label{gamma0}
\end{equation}
such that
\begin{equation}
\Gamma_n
\subset \Gamma\subset \Gamma_k, \qquad\mbox{for some}\ k>\frac {n+1}2,
\label{gamma1}
\end{equation}
and
\begin{equation}
\Gamma\ \mbox{is symmetric in the}\ \lambda_i.
\label{gamma1new}
\end{equation}
Let
\begin{equation}
f\in C^2(\Gamma)\cap
C^0(\overline \Gamma)\
\mbox{be  concave and symmetric in the}\ \lambda_i.
\label{hypo0}
\end{equation}
In addition, we assume that
\begin{equation}
f=0\ \mbox{on}\ \partial \Gamma;\quad
f_{\lambda_i}>0\ \mbox{on}\ \Gamma\ \forall\  1\le i\le n,
\label{aa5}
\end{equation}
and
\begin{equation}
\lim_{s\to\infty}f(s\lambda)=\infty,\qquad
\forall\ \lambda\in \Gamma.
\label{aa6}
\end{equation}

By (\ref{aa5}) and (\ref{aa6}), there exists a unique $
\bar b>0$ such that
\begin{equation}
f(\bar be)=1,
\label{b}
\end{equation}
where $e=(1,\cdots, 1)$.

\begin{cor}  For $n\ge 3$, let  $(f,\Gamma)$ satisfy
(\ref{gamma0}), (\ref{gamma1}), (\ref{gamma1new}),
 (\ref{hypo0}), (\ref{aa5})
and (\ref{aa6}), and let $u\in C^4(\Bbb R^n)$ be  a positive solution  of
\[
f(\lambda(A^u))=1,\quad \lambda(A^u)\in \Gamma, \quad\mbox{on}~\Bbb{R}^n.
\]
Then for some $\bar{x}\in\Bbb{R}^n$,
 and  some positive constant $a$,
$$
u(x)\equiv
 \left(\frac a {1+\frac 12 a^2\bar b|x-\bar{x}|^2}\right)^{\frac{n-2}{2}},
\qquad \forall\ x\in \Bbb R^n.
$$
\label{cor1}
\end{cor}

\noindent {\bf Proof of Theorem~\ref{highk}.}\ 
Since (\ref{Gammak}) implies the
superharmonicity of the positive function $u$ on $\Bbb{R}^n$, we have
$\liminf\limits_{|x|\to \infty}|x|^{n-2}u(x)>0$. Let
$w(x)=\frac{1}{|x|^{n-2}}u(\frac{x}{|x|^2})$ for
$x\in\Bbb{R}^n\setminus\{0\}$. Then $w$ 
is regular at $\infty$, $\liminf\limits_{|x|\to
0}w(x)>0$, and  $w$ satisfies
\[
F(A^w)=1,~~A^w\in U,\qquad\mbox{on}~\Bbb{R}^n\setminus\{0\}.
\]

Let
$\xi(x)=\frac{n-2}{2}w(x)^{-\frac{2}{n-2}}$.
Then, for some positive constant $C_1$,
\begin{equation}
0<\xi<C_1\quad\mbox{on}~B_2\setminus\{0\}.
\label{C0}
\end{equation}
By (\ref{Gammak}) and  lemma~6.3 in \cite{LL}, 
$\lambda(D^2\xi(x))\in\Gamma_k$ for $x\in
B_2\setminus\{0\}$. 
Let $P$ be any  hyperplane  which intersects $B_1$  but  does not pass
through the origin, and let $\xi_P$ be the restriction of
$\xi$ on $P$.   Then 
\[
\lambda(D^2\xi_P)\in\Gamma_{k-1}\subset\Bbb{R}^{n-1},\qquad\mbox{on}~P\cap
B_2,
\]
where $D^2\xi_P$ denotes $(n-1)\times (n-1)$ Hessian
of $\xi_P$, and $\lambda(D^2\xi_P)$ denotes
the eigenvalues of $D^2\xi_P$.
Here we have used the following property of $\Gamma_k$:
If $\lambda(M)\in \Gamma_k\subset \Bbb R^n$, then $\lambda(\hat M)
\in \Gamma_{k-1}\subset \Bbb R^{n-1}$ where
$\hat M_{ij}=M_{ij}$ for $1\le i,j\le n-1$.
  Since $k>\frac {n+1}2$, we have
$k-1>\frac{n-1}{2}$.  As in \cite{LL}, by using  theorem 2.7 in \cite{TW},
we have, for some  constants 
$\alpha \in (0,1)$ (depending only on $n$ and $k$) and
$C>0$  (depending only on $n$, $k$ and $C_1$), that
\begin{equation}
\|\xi\|_{C^{\alpha}(P\cap B_1)}\le C.
\label{Calpha0}
\end{equation}

For any $x, y\in B_1\setminus\{0\}$, we pick $z_i\in \Bbb R^n$ such that
$z_i\to 0$ and the line going through $x$ and $y+z_i$ does not
go through the origin.  Then $x$ and $y+z_i$ lies
on some hyperplane $P_i$ which does not go through the origin.
Thus, by (\ref{Calpha0}), 
$$
|\xi(x)-\xi(y+z_i)|\le C|x-(y+z_i)|^{\alpha}
$$
for some constant $C$ depending  only  on $n$, $k$ and $C_1$.
Sending $i$ to infinity, we have
$$
|\xi(x)-\xi(y)|\le C|x-y)|^{\alpha}.
$$
Therefore $\xi$ can be extended to a function
in $C^{\alpha}(B_1)$.

We distinguish into two cases.

\noindent {\bf Case 1.}\ $\xi(0)=0$.

\noindent {\bf Case 2.}\ $\xi(0)>0$.

In Case 1, $\lim_{ |x|\to \infty }(|x|^{n-2}u(x))=\infty$.
For every $x\in \Bbb R^n$, as in the proof of
lemma 2.1 in   \cite{LZ}, there exists
$\lambda_0(x)>0$ such that
$$
u_{x,\lambda}(y):=(\frac \lambda{|y-x|})^{n-2}
u(x+\frac {\lambda^2 (y-x)}{ |y-x|^2})
\le u(y), \ \forall\ 0<\lambda<\lambda_0(x), |y-x|\ge \lambda.
$$
Set, for $x\in \Bbb R^n$,
$$
\bar\lambda(x)=\sup\{ \mu>0\ |\
u_{x,\lambda}(y)\le u(y), \ \mbox{for all}\
|y-x|\ge \lambda, 0<\lambda \le \mu\}.
$$

\begin{lem} $\bar\lambda (x)=\infty$ for
all $x\in \Bbb R^n$.
\label{claim1}
\end{lem}

\noindent {\bf Proof of Lemma \ref{claim1}.}\
If $\bar\lambda (\bar x)<\infty$ for some $\bar x\in \Bbb R^n$.
Making a translation, we
may assume without loss of generality that
$\bar x=0$, and we still have
$\lim_{|x|\to\infty} (|x|^{n-2}u(x))=\infty$.
Thus, 
 there exists some
$R>\bar \lambda+9$ (we use notation
$\bar \lambda=\bar\lambda(0)$ such that
\begin{equation}
u_{\lambda}(y)< u(y),\qquad
\forall\ 0<\lambda\le \bar\lambda+2,
\ |y|\ge R,
\label{guarentee}
\end{equation}
where we have used notation $u_{\lambda}=u_{0,\lambda}$.

By the definition of $\bar\lambda$,
$$
u_{\bar \lambda}(y)\le u(y),\qquad \forall\
 |y|\ge \bar\lambda.
$$
Let $w_t:=tu+(1-t)u_{\bar \lambda}$, $0\le t\le 1$.
Then, as in the proof of lemma 2.1 in \cite{LL},
$$
L(u-u_{\bar \lambda})=0, \qquad \mbox{in}\
\Bbb R^n\setminus B_{ \bar\lambda },
$$
where
$$
L=a_{ij}(y)\partial_{ij}+b_i(y)\partial_i+c(y),
$$
$$
a_{ij}=-\frac 2{n-2}\int_0^1 w_t^{ -\frac {n+2}{n-2} }F_{ij}(A^{ w_t})dt,
$$
and $b_i$ and $c$ are continuous functions.

Using the Hopf Lemma and the
strong maximum principle as 
 in the proof of lemma 2.1 in \cite{LL},
we have
$$
(u-u_{\bar \lambda})(y)>0, \qquad \mbox{in}\
\Bbb R^n\setminus \overline B_{ \bar\lambda },
$$
and
$$
\frac { \partial (u-u_{\bar\lambda}) }
{\partial r}\bigg|_{ \partial B_{ \bar\lambda  } }
>0,
$$
where $\frac {\partial }{\partial r}$ denotes the 
outer normal differentiation.

The following argument is similar to the one used
in 
 the proof of lemma 2.2 in \cite{LZ}.
Since $ \partial B_{ \bar\lambda  }$
is compact,
$
\frac { \partial (u-u_{\bar\lambda}) }
{\partial r}\bigg|_{ \partial B_{ \bar\lambda  } }$ has a positive lower bound. 
Using the 
 $C^1$ regularity of $u$,
we can find  some $0<\delta<1$ such that
$$
\frac { \partial (u-u_{\lambda}) }
{\partial r}(y)>0,\qquad
\forall \bar\lambda\le \lambda\le \bar\lambda+\delta,
\lambda\le |y|\le \lambda+\delta.
$$
Since $(u-u_{\lambda})(y)=0 $ for $|y|=\lambda$, the above
implies 
$$
u_\lambda(y)\le u(y), \qquad 
\forall \bar\lambda\le \lambda\le \bar\lambda+\delta,
\lambda\le |y|\le \lambda+\delta.
$$

Since $(u_{\bar \lambda}-u)(y)<0$
for $\bar \lambda+\delta\le |y|\le R$, and since the set is compact,
there exists $\epsilon\in (0,\delta)$ 
such that
$$
u_{\lambda}(y)<u(y),\qquad
\forall\ \bar \lambda
\le \lambda\le \bar \lambda+\epsilon,
\bar \lambda+\delta\le |y|\le R.
$$
Here we have used the 
the continuity of $u$.

We have proved,
for the  $\epsilon$ above, that
$$
u_{\lambda}(y)\le u(y),
\qquad\forall\ \bar\lambda\le \lambda\le \bar\lambda
+\epsilon, \ |y|\ge \lambda.
$$
This violates the definition of $\bar\lambda$.  Lemma
\ref{claim1} is 
established.

\vskip 5pt
\hfill $\Box$
\vskip 5pt

It follows from  Lemma \ref{claim1} that
$$
u_{x,\lambda}(y)\le u(y),\qquad
\forall\ x\in \Bbb R^n, 0<\lambda<\infty, |y-x|\ge \lambda.
$$
This, together with some calculus lemma (see, e.g.,
lemma 11.2 in \cite{LZ}), implies that $u$ is a constant on
$\Bbb R^n$, thus $A^u\equiv 0$.  This is impossible because
of (\ref{F3}).  We have ruled out Case 1.

In Case 2, there exists some  constant $0<\delta<\frac 1 {20}$   
 such that
\begin{equation}
\delta\le w\le\frac{1}{\delta},\qquad\mbox{on}~B_{10\delta}.
\label{delta}
\end{equation}   

\begin{lem}
$$
\limsup_{|x|\to 0}
(|x||\nabla w(x)|)<\infty.
$$
\label{order1}
\end{lem}

\noindent {\bf Proof of
Lemma \ref{order1}.}\
For any  $0<r<5\delta$, let
$v(y):=w(ry)$ for $0<|y|<2$. Then $v$ satisfies
\begin{equation}
F(r^{-2} A^v)=1,\quad A^v\in U,\qquad\mbox{on}~B_2\setminus\{0\}.
\label{ev}
\end{equation}

For any $x\in B_{\frac 32}\setminus B_{\frac 34}$, as in the
proof of lemma 2.1 in \cite{LZ}, there exists $\lambda_0(x)\in (0, \frac 15)$
 such that
$$
v_{x,\lambda}(y)
:=(\frac{\lambda}{|y-x|})^{n-2}v(x+\frac{\lambda^2(y-x)}{|y-x|^2})
\le v(y),
 \ \ \forall\
y\in (B_2\setminus B_{\frac
12})\setminus B_{ \lambda }(x),
\ 0<\lambda\le \lambda_0(x).
$$

Set, for $x\in  B_{\frac 32}\setminus B_{\frac 34}$,
$$
\bar\lambda(x)=\sup\{ \mu>0\ |\
v_{x,\lambda}(y)\le v(y), \ \forall\ 
y\in \left(B_2\setminus B_{\frac
12}\right)\setminus B_{ \lambda }(x),
\ 0<\lambda\le \mu\}.
$$
Using the Hopf Lemma and the strong maximum principle as
in the proof of lemma 2.1 in \cite{LL}, and using 
the argument in Lemma \ref{claim1}, we know that  
 stopping at $\bar \lambda(x)$ is due to a boundary
touching, i.e.,
there exists
some  $y_0\in\partial
(B_2\setminus B_{\frac 12})$ such that
 $v_{x,\bar\lambda(x)}(y_0)=v(y_0)$,
i.e.
\[
(\frac{\bar \lambda(x)}{|y_0-x|})^{n-2}
w(rx+\frac{\lambda^2r(y_0-x)}{|y_0-x|^2})=w(r
y_0),
\] 
from which we deduce,
using (\ref{delta}),  that 
\[
\bar \lambda(x)^{n-2}= |y_0-x|^{n-2}\frac{w(ry_0)}{w(rx+\frac{\lambda^2
r(y_0-x)}{|y_0-x|^2})}\ge \delta^2 |y_0-x|^{n-2}\ge
4^{2-n}\delta^2.
\]
Thus we have shown that for any 
$x\in  B_{\frac 32}\setminus B_{\frac 34}$  and any
$0<\lambda<\frac 14
\delta^{\frac{2}{n-2}}$ we have 
\[
v_{x,\lambda}(y)\le v(y),\quad\forall~
y\in B_2\setminus B_{\frac
12},~|y-x|\ge \lambda.
\]
This and  some 
calculus lemma (see lemma 1 in \cite{LL2}) imply,
for some constant $C$ depending only on $\delta$, that
$$
|\nabla v(y)|\le C v(y)\qquad\forall\ |y|=1,
$$
i.e.,
\[
|\nabla w(ry)|\le C\frac{w(ry)}{r},\qquad\forall |y|=1.
\]
Since this holds for all $0<r<5\delta$, 
we have
\[
|z||\nabla w(z)|\le C w(z),\qquad\forall~0<|z|<5\delta.
\]
Lemma 
\ref{order1} is established.

\vskip 5pt
\hfill $\Box$
\vskip 5pt

Our next lemma provides estimates of the second derivatives of
$w$ near the origin.

\begin{lem} 
$$
\limsup_{|x|\to 0} (|x|^2 |\nabla^2 w(x)|)<\infty.
$$
\label{second}
\end{lem}

\noindent {\bf Proof of Lemma \ref{second}.}\
Let $\delta$ be as in the proof of Lemma \ref{order1}, $0<r<5\delta$,
and $v(y):=w(ry)$.
Then $v$ satisfies (\ref{ev}), i.e., 
$$
\tilde F(A^v)=r^2,\quad A^v\in \tilde U, \qquad
\mbox{on}\ B_2\setminus\{0\},
$$
where $\tilde U:=r^{2}U$ and
$\tilde F(M):=r^2F(r^{-2}M), M\in \tilde U$.
Clearly, 
  $(\tilde F,\tilde U)$ satisfies
(\ref{U1}), (\ref{U2}), (\ref{Gammak}),
(\ref{F1}),  (\ref{F2}), (\ref{F4})
 (with $(F,U)$ replaced by $(\tilde F, \tilde U)$) , and
$$
\sum_{i,j=1}^n
\tilde F_{ij}(M)M_{ij}\le \frac 1\gamma |M|^{1-\gamma} 
\sum_{i=1}^n \tilde F_{ii}(M),
\qquad\forall\
M\in \tilde U, \tilde F(M)=r^2,
|M|\ge 1.
$$

We know from  (\ref{delta}) and
Lemma \ref{order1} that
\[
v+|\nabla v|\le C\qquad\mbox{on}\ B_{\frac 32}\setminus B_{\frac 34} 
\]
for some constant $C$ independent of $r$.

Following, with
minor modification,  the computation in the proof of theorem 1.6 in 
\cite{LL} (with $F$ there replaced by our $\tilde F$, 
$v$ there replaced  by $-\frac 2{n-2}\log v$ with our $v$,
and keep in mind that $h$ there is a constant $r^2$;
for some earlier works 
on second derivative estimates, see remark 1.13 in \cite{LL}), we obtain
$$
|\nabla^2 v|\le C \quad\mbox{on}\ \partial B_1
$$
for some constant $C$ independent of $r$.
Lemma \ref{second}
follows immediately.

\vskip 5pt
\hfill $\Box$
\vskip 5pt

Since  $w\in C^\alpha(B_1)$
and since we have proved that
$$
\limsup_{|x|\to 0}(|x||\nabla w(x)|+|x|^2|\nabla ^2 w(x)|)<\infty,
$$
we can apply lemma  6.4 in \cite{LL} to obtain
$\limsup_{|x|\to 0}(|x|^{ 1-\frac \alpha 2}|\nabla w(x)|)<\infty$.
In particular,
$$
\lim_{|x|\to 0}( |x||\nabla w(x)| )=0.
$$

Now we are in a position to apply theorem 1.2 in \cite{LL}
(with $u_{0,1}$ there being our $w$) to conclude that
$u$ must be of the form (\ref{form}).
Theorem \ref{highk} is established.

\vskip 5pt
\hfill $\Box$
\vskip 5pt

\noindent{\bf Proof of Corollary \ref{cor1}.}\
Let
$$
U:=\{ M\in {\cal S}^{n\times n}\ |\
\lambda(M)\in \Gamma\},
$$
and
$$
F(M):=f(\lambda(M)),\qquad M\in U.
$$
To establish Corollary \ref{cor1}, 
we only need to verify that $(F,U)$ satisfies the hypothesis of
Theorem \ref{highk}.  
 These are well known to people in the field,
but for  convenience of the reader, we
provide some details. Since $\Gamma$ is an open subset of $\Bbb R^n$,
 $U$ is an
open subset of ${\cal S}^{n\times n}$.  Since orthogonal conjugation
does not change the set of eigenvalues and since $\Gamma$ is
symmetric in the $\lambda_i$, we know that $U$ satisfies (\ref{U1})
and $F$ satisfies (\ref{F1}).
Since $\Gamma_n\subset \Gamma$ and $\Gamma$ satisfies
(\ref{gamma0}), we know that $\lambda+\mu
=2(\frac {\lambda+\mu}2)\in \Gamma$ for all $\lambda\in \Gamma$ and
$\mu\in \overline \Gamma_n$.  For $M\in U$ and 
$N\in {\cal S}^{n\times n}$, let $\lambda_n(M)\ge \cdots\ge \lambda_1(M)$
 denote the eigenvalues of $M$, we know that
$$
\lambda_i(M)=\inf_{ dim\ K=i}
\sup_{ x\in X, \|x\|=1}
(x'Mx),\qquad 1\le i\le n.
$$
Similar formula holds for $M+N$.
Thus $\lambda_i(M+N)\ge \lambda_i(M)$
for all $1\le i\le n$.  Write $\lambda=(\lambda_1(M),
\cdots, \lambda_n(M))$ and $\mu =(\lambda_1(M+N)-\lambda_1(M),
\cdots, \lambda_n(M+N)-\lambda_n(M))$, then $\lambda\in \Gamma$ and
$\mu\in \overline \Gamma_n$, thus $\lambda+\mu
=(\lambda_1(M+N),
\cdots, \lambda_n(M+N))\in \Gamma$, i.e. $M+N\in U$.
So $U$ satisfies (\ref{U2}).
Since $\Gamma_U=\Gamma$, (\ref{Gammak}) follows from
(\ref{gamma1}).  Clearly, (\ref{F3}) follows from $f(0)=0$.
Property (\ref{F2}) and (\ref{F4}) can be deduced 
from the concavity of $f$ in $\Gamma$ and the fact
that $f_{\lambda_i}>0$ in $\Gamma$ for every
$1\le i\le n$, see e.g., \cite{CNS}.
For all $\lambda\in\Gamma$ satisfying
$f(\lambda)=1$, we have, using  the concavity of $f$ in $\Gamma$ and
the convexity of $\Gamma$, 
$$
1=f(\bar b e)\le f(\lambda)+\sum_{i}f_{\lambda_i}(\lambda)
(\bar b-\lambda_i)=1+\sum_{i}f_{\lambda_i}(\lambda)
(\bar b-\lambda_i),
$$
i.e.,
$$
\sum_{i}f_{\lambda_i}(\lambda)\lambda_i\le \bar b
\sum_{i}f_{\lambda_i}(\lambda).
$$
This, after diagonalizing $M$ by an orthogonal conjugation,
implies (\ref{F5}).
Corollary \ref{cor1} is established.

\vskip 5pt
\hfill $\Box$

\end{document}